\documentclass[draft,english]{article}
\usepackage{amsfonts,amsmath,amsthm,amscd,amssymb,latexsym}

    \usepackage[T2A]{fontenc}
    \usepackage[cp1251]{inputenc}
    \usepackage[russian]{babel}

\def\adb{\allowdisplaybreaks}
\def\a1{\alpha}
\def\e1{\varepsilon}
\def\de1{\delta}

\def\s1{\sigma}
\def\f1{\varphi}
\def\th{\theta}
\def\imi{\int\limits_0^t}

\def\imx{\int\limits_0^h}
\def\kt{\sqrt{\th(t)-\th(\tau)}}
\def\tat{\th(t)-\th(\tau)}
\def\m1{\sqrt{t^{\beta+1}-\tau^{\beta+1}}}
\def\bb{\beta+1}

\newcommand{\nn}{\nonumber}
 \newcommand{\D}{\displaystyle}
\newcommand{\lb}{\linebreak}

\begin{document}
{\bf An inverse problem for strongly degenerate heat equation}

\bigskip

{\bf Mykola Ivanchov, Nataliya Saldina}

\bigskip

Department of Mechanics and Mathematics, Lviv National University,
79002 Lviv, Ukraine
\bigskip
E-mail: ivanchov@franko.lviv.ua

 \bigskip
AMS Mathematics Subject Classification 35R30, 35K65
 \bigskip

  \bigskip
{\bf Abstract.} In this paper we consider an inverse problem for
determining time - dependent heat conduction coefficient which
vanishes at initial moment as a power $ t^{\beta}. $ The case of
strong degeneration ($ \beta \ge1$) is studied.
 To prove the existence of solution we employ the Schauder fixed point theorem.
The uniqueness of the solution is established too.
  \bigskip


{\bf 1. Introduction}
\bigskip

Degenerate parabolic problems arise in a lot of fields of natural
and social
 sciences and technology (see, e.g., [1-5]).  These  problems may  be divided
 on different classes accordingly to the way of degeneration  with
 respect either to spatial variables or to time variable, weak or strong
 degeneration.  Direct problems for degenerate parabolic  equations are
 sufficiently well studied. As examples, we can mention the works [6-12]. On
 the  other hand, inverse problems for non-degenerate parabolic equations are
  no less investigated [13-17]. However, inverse problems for degenerate
  partial differential equations are almost not considered. There are some
  works [18-20] dedicated to inverse problems for partial differential equations
  degenerating with respect to a spatial variable.

In this paper we consider an inverse problem for the heat equation
with unknown
 heat conduction coefficient depending on time variable $ t.$  It is supposed that the unknown
coefficient vanishes at the initial moment as a power  $ t^{\beta}.
$ The case of weak degeneration ($\beta< 1$) was studied in [21].
Here we investigate the case of strong degeneration ($\beta\ge 1$).

In a domain $ Q_T\equiv\{(x,t):0<x<h,0<t<T\} $ we consider the
following heat equation
\begin{equation}\label{1.1}
u_t=a(t)u_{xx}+f(x,t)
\end{equation}
with unknown coefficient  $ a(t)> 0,t\in(0,T], $ initial condition
\begin{equation}\label{1.2}
u(x,0)=\varphi (x),\quad x\in [0,h],
\end{equation}
boundary conditions
\begin{equation}\label{1.3}
u(0,t)=\mu_1(t), \quad  u(h,t)=\mu_2(t),\quad t\in [0,T],
\end{equation}
and overdetermination condition
\begin{equation}\label{1.4}
a(t)u_x(0,t)=\mu_3(t),\quad t\in [0,T].
\end{equation}

In this problem both $ u(x,t) $ and the coefficient $ a(t) $ are
unknown, are to be determined from data $ f(x,t), \f1(x),
\mu_1(t),\mu_2(t), \mu_3(t). $ As a solution of the problem
(\ref{1.1})-(\ref{1.4}) we mean a classical solution which is
defined as follows.

{\bf Definition. } {\it The pair of functions $ a(t) $ and $ u(x,t)
$ is a solution of (\ref{1.1})-(\ref{1.4}) if the following
conditions are fulfilled:
\begin{description}
\item[(a)] $ (a,u) \in  C[0,T]\times C^{2,1}( Q_T)\cap C(\overline
 Q_T),u_x(0,t)\in C(0,T];$
\item[(b)] $ a(t) $ is positive for $ t\in (0,T]; $
 \item[(c)] there exists the limit  $ \lim\limits_{t\to +0} \dfrac{a(t)}
 {t^{\beta }} >0 $ , $ \beta \ge1 $ - a given  number;
 \item[(d)] (\ref{1.1})-(\ref{1.4}) are satisfied.
\end{description}}

\smallskip

Note that the analogous problem for a non-degenerate heat equation
was for the first time studied in [22].

We establish the conditions of existence and uniqueness of solution
for the problem (\ref{1.1})-(\ref{1.4}) which are formulated in the
following theorem.

\smallskip

{\bf Theorem.} {\it Suppose that the following conditions hold:
\begin{description}
\item[1)] $ \varphi \in C^2[0,h];\mu_i \in C^1[0,T],i=1,2;\mu_3 \in C[0,T];
 f\in C^{2,0}(\overline Q_T);$
\item[2)]$ \varphi'(x)\ge 0,x\in [0,h];f(0,t)-\mu_1'(t)> 0,\mu_2'(t)-f(h,t)
\ge 0,t\in [0,T];\mu_3(t)> 0,
 t\in (0,T],$ the limit  $ \lim\limits_{t\to +0}
 \D\frac{\mu_3(t)}{t^{\frac{\beta+1}2}}>0 $ exists; $ f_x(x,t)\ge 0,(x,t)\in\overline Q_T;$
 \item[3)] $ \varphi(0)=\mu_1(0),\varphi(h)=\mu_2(0).$
\end{description}

Then there exists an unique solution of problem
(\ref{1.1})-(\ref{1.4}).}

 \smallskip

To prove the existence of solution, the Schauder fixed point theorem
is applied. The proof of the uniqueness of solution is divided in
two parts: first we establish it for a small time interval and after
this we prove a global (in time) uniqueness of solution of the
problem (\ref{1.1})-(\ref{1.4}).

The part of the paper that follows is composed of four sections. In
Section 2 the inverse problem (\ref{1.1})-(\ref{1.4}) is reduced to
an integral equation with respect to unknown coefficient $a(t).$ In
Section 3 we study the behavior of $ u_x(0,t) $ as  $t\to 0 $  and
we show that under the assumptions of the theorem the solution of
the integral equation is bounded from below and above
 by a power $t^{\beta} $ with coefficients which depend on given data. In Section 4
 we apply the Schauder fixed point theorem to the integral equation and we
 complete the proof of the existence of solution of the problem
 (\ref{1.1})-(\ref{1.4}).
In the first part of Section 5 we show that the integral equation
with respect to unknown coefficient $ a(t) $ admits at most one
solution on the interval $[0,\tilde t] $ where $ \tilde t>0 $ is, in
general, a small number defined by given data. Then using this
statement, we prove the uniqueness of solution of the problem
(\ref{1.1})-(\ref{1.4}) in whole.

\bigskip

 {\bf 2. Reduction of the problem (\ref{1.1})-(\ref{1.4}) to
an integral equation}

\bigskip

 Apply the overdetermination condition (\ref{1.4}) to obtain an equation
for the function  $ a(t)$. Denote by $ G_k(x,t,\xi,\tau), k=1,2, $
the Green functions of the first $ (k=1) $ and the second $ (k=2) $
boundary value problems for equation (\ref{1.1})
\begin{align}\label{2.1}
\adb
&G_k(x,t,\xi,\tau)=\frac 1{2\sqrt
{\pi(\theta(t)-\theta(\tau))}}\sum\limits_{n= -\infty
}^{\infty}\biggl(\exp\biggl(-\frac{(x-\xi+2nh)^2}{4(
\theta(t)-\theta(\tau))}\biggr)\nn\\
&+(-1)^k\exp\biggl(-\frac{(x+\xi+2nh)^2}{4(
\theta(t)-\theta(\tau))}\biggr)\biggr),
\end{align}
 where $
\theta(t)=\D\int\limits_0^ta(\tau )d\tau.$
 Temporarily assuming that the function $ a(t) $ is known, we can write the
 solution of direct problem  (\ref{1.1})-(\ref{1.3}) with the aid of the Green function
\begin{align}\label{2.2}
\adb
&u(x,t)=\int\limits_0^hG_1(x,t,\xi,0)\varphi(\xi)d\xi+\int\limits_0^tG_{1\xi}
(x,t,0,\tau)a(\tau)\mu_1(\tau)d\tau\nn\\
&-\int\limits_0^tG_{1\xi}(x,t,h,\tau)a(\tau)\mu_2(\tau)d\tau+\int\limits_0^t
\int\limits_0^hG_1(x,t,\xi,\tau)f(\xi,\tau)d\xi d\tau.
\end{align}

Evaluate the first derivative $ u_x(x,t) $, taking into account the
relationships
\begin{equation}\label{2.3}
G_{1x}(x,t,\xi,\tau)=-G_{2\xi}(x,t,\xi,\tau)\quad\text{and}\quad
G_{2\xi\xi}=-\dfrac {G_{2\tau}(x,t,\xi,\tau)}{a(\tau)},
\end{equation}
 and
integrating by parts with using compatibility condition. We obtain
\begin{align}\label{2.4}
\adb
&u_x(x,t)=\int\limits_0^hG_2(x,t,\xi,0)\varphi'(\xi)d\xi+\int\limits_0^tG_2
(x,t,0,\tau)(f(0,\tau)-\mu_1'(\tau))d\tau+\nn\\
&+\int\limits_0^tG_2(x,t,h,\tau)(\mu_2'(\tau)-f(h,\tau))d\tau+\int\limits_0^t
\int\limits_0^hG_2(x,t,\xi,\tau)f_{\xi}(\xi,\tau)d\xi d\tau.
\end{align}

We substitute this expression into overdetermination condition
(\ref{1.4}) and we come to the equation for $ a(t) $ :
\begin{align}\label{2.5}
\adb
&a(t)=\mu_3(t)\biggl(\int\limits_0^hG_2(0,t,\xi,0)\varphi'(\xi)d\xi+\int\limits_0^tG_2
(0,t,0,\tau)(f(0,\tau)-\mu_1'(\tau))d\tau\nn\\
&+\int\limits_0^tG_2(0,t,h,\tau)(\mu_2'(\tau)-f(h,\tau))d\tau+\int\limits_0^t
\int\limits_0^hG_2(0,t,\xi,\tau)f_{\xi}(\xi,\tau)d\xi d\tau
\biggr)^{-1}, t\in [0,T].
\end{align}

Taking into account the conditions of the theorem it is easy to
verify that the function $ a(t) $ is positive on $ (0,T] $ and
belongs to $ C(0,T] $.

\bigskip
{\bf 3. A priori estimates}
\bigskip

In order to prove the existence of the solution of equation
(\ref{2.5}) we apply the Schauder fixed point theorem. First of all,
we estimate the solution of equation (\ref{2.5}). As a consequence
of the second condition of the theorem and the explicit
representation of the function  $ G_2(x,t,\xi,\tau) $ we obtain the
inequality
\begin{equation}\label{3.1}
u_x(0,t) \ge \imi\frac
{f(0,\tau)-\mu_1'(\tau)}{\sqrt{\pi(\th(t)-\th(\tau))}}d\tau.
\end{equation}
This allows us to write the following inequality for the function $
a(t) $
\begin{equation}\label{3.2}
a(t)\le \frac{\mu_3(t)}{\D\imi\dfrac
{f(0,\tau)-\mu_1'(\tau)}{\sqrt{\pi(\th(t)-\th(\tau))}}d\tau} .
\end{equation}
Let
\begin{equation}\label{3.3}
 a_0(t)\equiv\D\frac{a(t)}{t^{\beta}},\quad  a_{\max}(t)\equiv\max\limits_{0\le \tau\le
  t}a_0(\tau).
 \end{equation}
 Then from (\ref{3.2}) we find
$$
a_0(t)\le\frac{\sqrt{\pi}\mu_3(t)}{\sqrt{\beta+1}t^{\beta}\D\imi\dfrac
{f(0,\tau)-\mu_1'(\tau)}{\m1}d\tau}\sqrt{a_{\max}(t)} .
$$
Denote
\begin{equation}\label{3.4}
H(t)\equiv\frac{\sqrt{\pi}\mu_3(t)}{\sqrt{\beta+1}t^{\beta}\D\imi\dfrac
{f(0,\tau)-\mu_1'(\tau)}{\m1}d\tau}.
\end{equation}
 It follows from the
conditions of the theorem that the function $ H(t) $ is positive on
$ (0,T] $ and belongs to $ C(0,T] $. Establish the existence of the
limit $ \lim\limits _{t\to +0}H(t).$ To this end, we apply the
theorem on average and change of variables $ z=\dfrac {\tau}{t} $:
\begin{align*}
\adb & \lim\limits_{t\to +0}H(t)=\lim\limits_{t\to
+0}\frac{\sqrt{\pi}\mu_3(t)} {\sqrt{\beta+1}t^{\beta}(f(0,\overline
t)-\mu_1'(\overline t))\D\imi
\dfrac{d\tau}{\m1}}=\\
&=\sqrt{\frac{\pi}{\beta+1}}\lim\limits_{t\to +0}\frac{\mu_3(t)}
{t^{(\beta+1)/2}(f(0,\overline t)-\mu_1'(\overline
t))\D\int\limits_0^1 \dfrac{dz}{\sqrt{1-z^{\beta+1}}}},\quad
\overline t \in[0,t].
\end{align*}
Denote
\begin{equation}\label{3.5}
\D\int\limits_0^1\frac{dz}{\sqrt{1-z^{\beta+1}}}=I_1,
\end{equation}
then we obtain $ \lim\limits_{t\to
+0}H(t)=\dfrac{\sqrt{\pi}M}{\sqrt{\beta+1} (f(0,0)-\mu_1'(0))I_1}>0,
$ where $ M=\lim\limits_{t\to +0}\dfrac{\mu_3(t)}
{t^{\frac{\beta+1}2}}.$ Applying (\ref{3.4}) we come to the
inequality $ a_0(t)\le H(t)\sqrt{a_{\max}(t)}, $
 which leads to an estimation of $ a_{\max}(t) $ from above
\begin{equation}\label{3.6}
 a_{\max}(t)\le H^2_{\max}(t)<\infty,\quad t\in[0,T],
\end{equation}
where $ H_{\max}(t)\equiv\max\limits_{0\le\tau\le t}H(\tau).$ Taking
into account the existence of the limit $ \lim\limits_{t\to +0}
\dfrac{\mu_3(t)} {t^{\frac{\beta+1}2}} $ and notation (\ref{3.5}),
we estimate $ H(t) $
 \begin{equation}\label{3.7}
  H(t)\le \frac{\sqrt{\pi}M_1}{\sqrt{\beta+1}\min\limits_{[0,T]}
 (f(0,t)-\mu_1'(t))
I_1}\equiv H_1,\quad\text{where}\quad
M_1=\max\limits_{[0,T]}\frac{\mu_3(t)} {t^{\frac{\beta+1}2}}.
\end{equation}
Then the estimate of  $ a(t) $ from above follows:
\begin{equation}\label{3.8}
a(t)\le H^2_{\max}(t)t^{\beta}\le H_1^2t^{\beta},\quad t\in[0,T].
\end{equation}

Estimate $ a(t) $ from below. To this end, provide some estimates
for expression in the denominator of (\ref{2.4}). From the equality
$ \D\int\limits_0^hG_2(x,t,\xi,\tau)d\xi=1 $ the first and the
fourth summands are estimated
$$
\int\limits_0^hG_2(0,t,\xi,\tau)\f1'(\xi)d\xi\le C_1,\quad
\int\limits_0^t\int
\limits_0^hG_2(0,t,\xi,\tau)f_{\xi}(\xi,\tau)d\xi d\tau \le C_2.
$$
It can be easily verified that
$$
\int\limits_0^tG_2(0,t,h,\tau)(\mu_2'(\tau)-f(h,\tau))d\tau\le C_3,
$$
where $ C_1,C_2,C_3>0 $ --- the constants determined by the problem
data. Transform the second summand, separating out of the series the
term that corresponds to $ n=0 : $
 \begin{align*} & \imi
G_2(0,t,0,\tau)(f(0,\tau)-\mu_1'(\tau))d\tau=
\frac1{\sqrt{\pi}}\imi\frac{f(0,\tau)-\mu_1'(\tau)}{\kt}d\tau\\
&+\frac2{\sqrt{\pi}}\imi\frac{f(0,\tau)-\mu_1'(\tau)}{\kt}\sum\limits_{n=1}^
{\infty}\text{exp}\left(-\frac{n^2h^2}{\tat}\right)d\tau.
\end{align*}
Since the integral function of the last summand has no singularities
and is bounded, we estimate it by a constant. Denote $
a_{\min}(t)\equiv\min\limits_{0\le \tau\le t}a_0(\tau). $ By formula
(\ref{2.4}) and previous estimates we obtain
\begin{align*} & \adb
a_0(t)\ge\frac{\mu_3(t)}{t^{\beta}\left(C_4+\dfrac{\sqrt{\beta+1}}{\sqrt{\pi
 a_{\min}(t)}}\D\imi\frac{f(0,\tau)-\mu_1'(\tau)}{\m1}d\tau\right)}\\
 &\ge\frac
 {\sqrt{a_{\min}(t)}}{\dfrac{C_5t^{\beta}}{\mu_3(t)}+\dfrac{t^{\beta}
 \sqrt{\beta+1}}{\sqrt{\pi}\mu_3(t)}\D\imi\frac{f(0,\tau)-\mu_1'(\tau)}{\m1}
 d\tau}.
\end{align*}
 Taking into account (\ref{3.4}), reduce the inequality to the form
$$
a_0(t)\ge\frac{\sqrt{a_{\min}(t)}}{\dfrac{C_5t^{\beta}}{\mu_3(t)}+\dfrac1{H(t)}}
=\frac{\sqrt{a_{\min}(t)}H(t)}{\dfrac{C_5t^{\beta}H(t)}{\mu_3(t)}+1}.
$$
From (\ref{3.7}) it follows
$$
\frac{C_5t^{\beta}H(t)}{\mu_3(t)}\le C_6t^{\frac{\beta-1}2}.
$$
Hence, we obtain
$$
a_0(t)\ge\frac{\sqrt{a_{\min}(t)}H(t)}{C_6t^{\frac{\beta-1}2}+1}.
$$
From here we establish the estimate
\begin{equation}\label{3.9}
a_{\min}(t)\ge\frac{H^2_{\min}(t)}{(C_6t^{\frac{\beta-1}2}+1)^2},
\end{equation}
where $ H_{\min}(t)\equiv\min\limits_{0\le\tau\le t}H(\tau).$
Finally, we have for $ a(t) $:
\begin{equation}\label{3.10}
0<A_0\le\frac{H^2_{\min}(t)}{(C_6t^{\frac{\beta-1}2}+1)^2}\le
\frac{a(t)}{t^{ \beta}}\le H^2_{\max}(t)\le A_1<\infty,\quad
t\in[0,T].
\end{equation}
 Therefore, we have established a priori estimates for
the solution (\ref{2.5}). Having a priori estimates of solution of
equation (\ref{2.5}), we can apply to it
 the Schauder fixed point theorem.

 \bigskip
{\bf 4. Existence of solution}
\bigskip

We consider the equation (\ref{2.5}) as an operator equation $
a(t)=Pa(t) $ with respect to $ a(t) $ and the operator $ P $ is
defined by equality $ Pa(t)=\D\frac {\mu_3(t)}{u_x(0,t)}. $ Denote $
\mathcal{ N}=\{a\in C[0,T]:A_0\le\dfrac{a(t)}{t^{\beta}} \le A_1\}.$
As a consequence of a priori estimates (\ref{3.10}) the operator $ P
$ maps
 $ \mathcal{N}  $ into $ \mathcal{N}  $. We are going to show that
  the set $ P\mathcal{N} $ is compact
or equivalently, by Arzela-Ascolli theorem, $ P\mathcal{N} $ is
uniformly bounded and equicontinuous. We have to establish that
$\forall\varepsilon>0 \ \exists \
 \delta>0 $ such that
\begin{equation}\label{4.1}
  |Pa(t_2)-Pa(t_1)|<\e1  \quad\text{for arbitrary} \quad  |t_2-t_1|<\delta,
  \quad a(t)\in \mathcal{N}.
 \end{equation}
 As $ Pa(t)= a(t) $ and $\dfrac{a(t)}{t^{\beta}}\le A_1 $  for all $ a(t) \in
 \mathcal{N} $ we conclude that for arbitrary $ \varepsilon>0$ there exists
 sufficiently small number $ t^*>0 $, such that the  inequality
  $$
  |Pa(t)|<\varepsilon,   \quad   0\le t\le t^*,
 $$
holds.

Establish the inequality (\ref{4.1}) in the case when $ t_i>t^*,
i=1,2.$ Assume $ t_2>t_1 $. Consider one of the summand which is
contained in (\ref{4.1}):
\begin{align*}
&R_1=\left|\int\limits_0^{t_2}G_2(0,t_2,0,\tau)(f(0,\tau)-\mu_1'(\tau))d\tau
-\int \limits_0^{t_1}G_2(0,t_1,0,\tau)(f(0,\tau)-\mu_1'(\tau))
d\tau\right|\\
&\le \left|\int \limits_0^{t_1}(G_2(0,t_2,0,\tau)-G_2(0,t_1,0,\tau))
(f(0,\tau)-\mu_1'(\tau))d\tau \right|\\
&+\left|\int\limits_{t_1}^{t_2}G_2(0,t_2,0,\tau)
(f(0,\tau)-\mu_1'(\tau))d\tau\right|\equiv R_{1,1}+R_{1,2}.
\end{align*}
The estimate of $ G_2(0,t,0,\tau) $ [17] allows us to write
\begin{align*}
&R_{1,2}\le\frac{\max\limits_{[0,T]}(f(0,t)-\mu_1'(t))}{\sqrt{\pi}}\int\limits_{t_1}^{t_2}\left(
 \frac1 {\sqrt{\th(t_2)-\th(\tau)}}+C_7\right)d\tau\\
&\le C_8\int\limits_{t_1}^{t_2}\frac
 {d\tau}{\sqrt {\th(t_2)-\th(\tau)}}+C_9(t_2-t_1).
 \end{align*}
Apply the definition of the set $\mathcal{N} $  for investigation of
the first summand in $ R_{1,2}:$
$$
\int\limits_{t_1}^{t_2}\frac{d\tau}{\sqrt {\th(t_2)-\th(\tau)}}\le
\sqrt{\frac{\beta+1}
{A_0}}\int\limits_{t_1}^{t_2}\frac{d\tau}{\sqrt{t_2^{\bb}-\tau^{\bb}}}\le\sqrt
{\frac{\beta+1}{A_0t_2^{\beta}}}
\int\limits_{t_1}^{t_2}\frac{d\tau}{\sqrt{t_2-\tau}}\le
C_{10}\sqrt{t_2-t_1}.
$$
Finally, we have
$$
  R_{1,2}\le C_{11}\sqrt{t_2-t_1}+C_9(t_2-t_1).
  $$
From $ R_{1,1} $, using the Green function representation and
separating out
 of the series the term that corresponds to $ n=0 $, we have
\begin{align*}
&R_{1,1}\le
\frac1{\sqrt{\pi}}\max\limits_{[0,T]}(f(0,t)-\mu_1'(t))\left(\int
\limits_0^{t_1}\left|\frac1{\sqrt{\th(t_2)-\th(\tau)}}-\frac1{\sqrt{\th(t_1)-
\th(\tau)}}\right| d\tau\right.\\
&+2\int\limits_0^{t_1}\left|\frac1{\sqrt{\th(t_2)-\th (\tau)}}
\sum\limits_{n=1}
^{\infty}\text{exp}\left(-\frac{n^2h^2}{\th(t_2)-\th(\tau)}\right)
\right.\\
&-\left.\left.\frac1{\sqrt{\th(t_1)
-\th(\tau)}}\sum\limits_{n=1}^{\infty}\text{exp}
\left(-\frac{n^2h^2}{\th(t_1)-\th(\tau)}\right)\right|d\tau
\right)\equiv
 R_{1,1,1}+R_{1,1,2}.
\end{align*}
Note that
\begin{align*}
\adb
&\frac1{\sqrt{\th(t_2)-\th(\tau)}}-\frac1{\sqrt{\th(t_1)-\th(\tau)}}=\\
&=\frac{\th(t_1)-\th(t_2)}{\sqrt{(\th(t_2)-\th(\tau))(\th(t_1)
-\th(\tau))}(\sqrt{\th(t_1)-\th(\tau)}+\sqrt{\th(t_2)-\th(\tau)})}.
\end{align*}
Then we employ the definition of the set $ \mathcal{N} $:
\begin{equation*}
R_{1,1,1}\le \frac1{\sqrt{\pi}}\max\limits_{[0,T]}(f(0,t)-
\mu_1'(t))\frac{A_1}{2A_0^{3/2}}\int\limits_0^{t_1}\left(\frac1{\sqrt{t_1^
{\bb}-\tau^{\bb}}}-\frac1{\sqrt{t_2^{\bb}-\tau^{\bb}}}\right)d\tau.
\end{equation*}
Transform this inequality by the change of variable $
z=\dfrac{\tau}{t_1} $:
$$
R_{1,1,1}\le\frac{C_{12}}{t_1^{\frac{\beta-1}2}}\left(\int\limits_0^1\frac{dz}
{\sqrt{1-z^{\bb}}}-\int\limits_0^1\frac{dz}
{\sqrt{\left({t_2}/{t_1}\right)^{\bb}-z^{\bb}}}\right).
$$
Denote
\begin{equation}\label{4.2}
 I(\omega)\equiv\D\int\limits_0^1\dfrac{dz}{\sqrt{\omega-z^{\bb}}},\quad
 \text {where}\quad \omega\in[1,\infty).
\end{equation}
Obviously $ \D\int\limits_0^1\dfrac{dz}
{\sqrt{\omega-z^{\bb}}}\le\D\int\limits_0^1\dfrac{dz}{\sqrt{1-z}}. $
It follows that $ I(\omega) $ is continuous on $ [1,\infty) $. This
means that for arbitrary $ \varepsilon>0 $ there exists $ \delta_1>0
$, such that $ R_{1,1,1} <\varepsilon \ \text{ for }\
|t_2-t_1|<\delta_1$.

We represent $ R_{1,1,2} $  in the form:
$$
R_{1,1,2}=\frac2{\sqrt{\pi}}\max\limits_{[0,T]}(f(0,t)-\mu_1'(t))\int\limits_0^{t_1}
\left|\int\limits_{\th(t_1)-\th(\tau)}^{\th(t_2)-\th(\tau)}\frac
d{dz}
\left(\frac1{\sqrt{z}}\sum\limits_{n=1}^{\infty}\exp\left(-\frac{n^2h^2}z
\right)\right)dz\right|d\tau.
$$
Estimating $ R_{1,1,2} $ by means of inequality $ x^{n}e^{-x^2}\le
M_n<\infty, x\in[0,\infty),n\in N, $ we come to the estimate
$$
R_{1,1,2}\le C_{13}|\th(t_2)-\th(t_1)|\le C_{14}|t_2^{\bb}-
t_1^{\bb}| <\varepsilon \quad\text{for}\quad|t_2-t_1|<\delta_2,
$$
where $ \delta_2>0 $ depends on $ \varepsilon$ and known data.
Others summands occuring in (\ref{4.1}) can be estimated
analogously. It follows that the conditions of the Schauder theorem
for equation (\ref{2.5}) hold, and, therefore there exists the
solution $ a=a(t) $ of the equation (\ref{2.5}), which belongs to $
C[0,T] $.

Substituting $ a(t) $ in (\ref{2.2}), we obtain a solution $
u=u(x,t) $ of direct problem (\ref{1.1})-(\ref{1.3}) which possesses
necessary smoothness.

\bigskip
{ \bf 5. Uniqueness of solution}
\bigskip

To prove the uniqueness of the solution of the problem
(\ref{1.1})-(\ref{1.4}), suppose that $ (a_i(t),\lb u_i(x,t)),i=1,2
$ are two solutions of the problem (\ref{1.1})-(\ref{1.4}). Denote $
a(t)\equiv{a_1(t)-a_2(t)}, \lb u(x,t)\equiv {u_1(x,t)- u_2(x,t)} $.
From the overdetermination condition (\ref{1.4}) we obtain
\begin{equation}\label{5.1}
a(t)=\frac{\mu_3(t)(u_{2x}(0,t)-u_{1x}(0,t))}{u_{2x}(0,t)u_{1x}(0,t)}.
\end{equation}
or, after using the notation (\ref{3.3}),
$$
a_0(t)=\frac{\mu_3(t)}{t^{\beta}u_{1x}(0,t)}\frac{\mu_3(t)}{t^{\beta}u_{2x}
(0,t)}\frac{t^{\beta}(u_{2x}(0,t)-u_{1x}(0,t))}{\mu_3(t)}.
$$
Next, we apply the equality $
\dfrac{\mu_3(t)}{t^{\beta}u_{ix}(0,t)}=\dfrac
{a_i(t)}{t^{\beta}},\quad i=1,2, $ and the estimate (\ref{3.6}):
\begin{equation}\label{5.2}
|a_0(t)|\le
H_{\max}^4(t)\frac{t^{\beta}|u_{2x}(0,t)-u_{1x}(0,t)|}{\mu_3(t)}.
\end{equation}
Estimate one of the summands which is contained in the expression $
|u_{2x}(0,t)-u_{1x}(0,t)|. $ Denote
\begin{align*}
\adb
&I_1\equiv\frac1{\sqrt{\pi}}\imi(f(0,\tau)-\mu_1'(\tau))\left(\frac1{\sqrt{\th_2(t)-
\th_2(\tau)}}-\frac1{\sqrt{\th_1(t)-\th_1(\tau)}}\right)d\tau\\
&+\frac2{\sqrt{\pi}}\imi(f(0,\tau)-\mu_1'(\tau))
\sum\limits_{n=1}^{\infty}
\left(\frac1{\sqrt{\th_2(t)-\th_2(\tau)}}\exp\left(-\frac{n^2h^2}{\th_2(t)-
\th_2(\tau)}\right)\right.\\
&-\left.\frac1{\sqrt{\th_1(t)-\th_1(\tau)}}\exp\left(-\frac{n^2h^2}{\th_1(t)-
\th_1(\tau)}\right)\right)d\tau\equiv I_{1,1}+I_{1,2}.
\end{align*}
Represent the second summand in the form:
$$
\adb
|I_{1,2}|\le\frac2{\sqrt{\pi}}\max\limits_{[0,T]}(f(0,t)-\mu_1'(t))\left|\imi
d\tau
\int\limits_{\th_1(t)-\th_1(\tau)}^{\th_2(t)-\th_2(\tau)}\left|\frac
d{dz}
\left(\frac1{\sqrt{z}}\sum\limits_{n=1}^{\infty}\exp\left(-\frac{n^2h^2}z\right)
\right)\right|dz\right|.
$$
From boundedness of the integrand  we conclude
\begin{align*} \adb
&|I_{1,2}|\le C_{15}\imi
|\th_2(t)-\th_2(\tau)-\th_1(t)+\th_1(\tau)|d\tau=
C_{15}\imi d\tau\int\limits_{\tau}^t|a_2(\sigma)-a_1(\sigma)|d\sigma=\\
&=C_{15}\imi d\tau\int\limits_{\tau}^t|a_0(\s1)|\s1^{\beta}d\s1\le
 C_{15}t^{\beta+2}\tilde
a_{\max}(t),
\end{align*}
where $ \tilde a_{\max}(t)\equiv\max\limits_{0\le \tau\le
t}|a_0(\tau)|.$

Transform $ I_{1,1} $ to the form
$$
\adb
I_{1,1}=\frac1{\sqrt{\pi}}\imi\frac{(f(0,\tau)-\mu_1'(\tau))
(\th_1(t)-\th_1
(\tau)-\th_2(t)+\th_2(\tau))}{\sqrt{(\th_2(t)-\th_2(\tau))
(\th_1(t)-\th_1(\tau))}
(\sqrt{\th_1(t)-\th_1(\tau)}+\sqrt{\th_2(t)-\th_2(\tau)})}d\tau.
$$
Taking into account (\ref{3.10}), we obtain
\begin{align*}
\adb &\left|\th_1(t)-\th_1(\tau)-\th_2(t)+\th_2(\tau)\right|\le
\int\limits_{\tau}^t |a_0(\s1)|\s1^{\beta}d\s1\le
\frac{t^{\bb}-\tau^{\bb}}{\bb}\tilde a_{\max}(t),\\
&\th_i(t)-\th_i(\tau)=
\int\limits_{\tau}^ta_{i0}(\s1)\s1^{\beta}d\s1\ge
\frac{H_{\min}^2(t)}{(C_6t^{\frac{\beta-1}2}+1)^2}\frac{t^{\bb}-\tau^{\bb}}
{\bb},\quad i=1,2,
\end{align*}
where $ a_{i0}=\D\frac{a_i(t)}{t^{\beta}}, i=1,2.$

Finally, we have the estimate of $ I_{1,1} $
$$
|I_{1,1}|\le\frac{\sqrt{\bb}(C_6t^{\frac{\beta-1}2}+1)^3}{2\sqrt{\pi}
H_{\min}^3(t)}\tilde a_{\max}(t)\imi
\frac{f(0,\tau)-\mu_1'(\tau)}{\m1}d\tau,
$$
or if we use (\ref{3.4})
$$
|I_{1,1}|\le\frac{(C_6t^{\frac{\beta-1}2}+1)^3}{2H_{\min}^4(t)}\frac{\mu_3(t)}
{t^{\beta}}\tilde a_{\max}(t).
$$
Other summands in the expression $ |u_{2x}(0,t)-u_{1x}(0,t)| $ are
estimated analogously to $ I_{1,2} $. We continue to estimate
(\ref{5.2}) as follows:
\begin{equation}\label{5.3}
 \tilde
a_{\max}(t)\le\frac{(C_6t^{\frac{\beta-1}2}+1)^3H_{\max}^4(t)}{2H_{\min}^4
(t)}\tilde a_{\max}(t)+S(t)\tilde a_{\max}(t),
\end{equation}
where $ S(t) $ is the sum of terms  depending on $ t $ which
vanishes for $ t=0 $. The existence of the limit $ \lim\limits_{t\to
+0}H(t)>0 $ implies the existence of the limit $
\lim\limits_{t\to+0}\dfrac{(C_6t^{\frac{\beta-1}2}+1)
^3H_{\max}^4(t)}{2H_{\min}^4(t)}=\dfrac12 $. Therefore, there exists
$ t_1:0<t_1\le T $ such that the following inequality holds:
\begin{equation}\label{5.4}
\frac{(C_6t^{\frac{\beta-1}2}+1)^3H_{\max}^4(t)}{2H_{\min}^4(t)}\le\frac34,
\quad t\in[0,t_1].
\end{equation}
Hence, we can rewrite (\ref{5.3}) in the form:
$$
 \tilde a_{\max}(t)\left(\frac14-S(t)\right)\le0,\quad t\in[0,t_1].
$$
There exists $ t_2:0<t_2\le T, $ such that $ \frac14-S(t)>0 $ for
any $ t\in[0,t_2]. $ We come to the contradiction. This implies that
 $ a_1(t) \equiv a_2(t) $ when $ t\in [0,\tilde t] $ with $ \tilde t=\min(t_1,t_2).$

Now we establish the uniqueness of the solution problem
(\ref{1.1})-(\ref{1.4}) for any $ t\in [0,T]. $ The functions $
(a(t), u(x,t)) $ satisfy the conditions
 \begin{align}
\adb
&u_t=a_1(t)u_{xx}+a(t)u_{2xx},\quad (x,t)\in Q_T,\label{5.5} \\
&u(x,0)=0,\quad x\in[0,h], \label{5.6}   \\
&u(0,t)=u(h,t)=0,\quad t\in[0,T], \label{5.7}\\
&a_1(t)u_x(0,t)=-a(t)u_{2x}(0,t),\quad t\in[0,T].\label{5.8}
\end{align}
Denote by $ G_1^{(i)}(x,t,\xi,\tau), i=1,2, $ the Green functions of
the equations $ u_t=a_i(t)u_{xx}$ with the boundary condition
(\ref{5.7}). The solution of the  (\ref{5.5})-(\ref{5.7}) can be
written with the aid of $ G_1^{(1)}(x, t,\xi,\tau) $
\begin{equation}\label{5.9}
u(x,t)=\imi\imx
G_1^{(1)}(x,t,\xi,\tau)a(\tau)u_{2\xi\xi}(\xi,\tau)d\xi d\tau.
\end{equation}
Substituting (\ref{5.9}) into (\ref{5.8}), we obtain the integral
equation for $ a(t) $
\begin{equation}\label{5.10}
a(t)u_{2x}(0,t)=-a_1(t)\imi\imx  G_{1x}^{(1)}(0,t,\xi,\tau)a(\tau)
u_{2\xi\xi}(\xi,\tau)d\xi d\tau,\quad t\in[0,T].
\end{equation}
Using the notation (\ref{3.3}), it is easy to see that
\begin{equation}\label{5.11}
a_0(t)=-\frac{a_1(t)}{t^{\beta}u_{2x}(0,t)}\imi\imx
G_{1x}^{(1)}(0,t,\xi,\tau)
a_0(\tau)\tau^{\beta}u_{2\xi\xi}(\xi,\tau)d\xi d\tau,
\end{equation}
or
\begin{equation}\label{5.12}
a_0(t)=\imi K(t,\tau)a_0(\tau)d\tau,\quad t\in[0,T],
\end{equation}
where
\begin{equation}\label{5.13}
K(t,\tau)\equiv-\frac{a_1(t)\tau^{\beta}}{t^{\beta}u_{2x}(0,t)}\imx
G_{1x}^{(1)} (0,t,\xi,\tau)u_{2\xi\xi}(\xi,\tau)d\xi.
\end{equation}
Establish that the kernel $ K(t,\tau) $ is integrable, using the
fact that $ u_2(x,t) $ is the solution of the problem
(\ref{1.1})-(\ref{1.4}). From (\ref{2.2}) we find the derivative
\begin{align}\label{5.14}
\adb
 &u_{2xx}(x,t)=\imx G_{1}^{(2)}(x,t,\xi,0)\f1''(\xi)d\xi+\imi G_{1\xi}^{(2)}
 (x,t,0,\tau)(\mu_1'(\tau)-f(0,\tau))d\tau\nonumber\\
 &+\imi G_{1\xi}^{(2)}(x,t,h,\tau)(f(h,\tau)-\mu_2'(\tau))d\tau+\imi\imx
 G_1^{(2)}(x,t,\xi,\tau)f_{\xi\xi}(\xi,\tau) d\xi d\tau .
\end{align}
 Substituting (\ref{5.14}) into (\ref{5.13}), we obtain:
\begin{align*}
&K(t,\tau)=-\frac{a_1(t)\tau^{\beta}}{t^{\beta}u_{2x}(0,t)}\imx
G_{1x}^{(1)}
 (0,t,\xi,\tau)\left(\imx G_1^{(2)}(\xi,\tau,\eta,0)\f1''(\eta)d\eta
 \right.\\
 &+\left.\int\limits
 _0^{\tau} G_{1\eta}^{(2)}(\xi,\tau,0,\s1)(\mu_1'(\s1)-f(0,\s1))d\s1+\int
 \limits_0^ {\tau} G_{1\eta}^{(2)}(\xi,\tau,h,\s1)(f(h,\s1)-\mu_2'(\s1))d\s1
\right.\\
&+\left.\int \limits_0^ {\tau}\imx
G_1^{(2)}(\xi,\tau,\eta,\s1)f_{\eta\eta}(
 \eta,\s1)
 d\eta d\s1\right)d\xi\equiv-\frac{a_1(t)\tau^{\beta}}{t^{\beta}u_{2x}(0,t)}
 \sum\limits_{i=1}^4K_i(t,\tau).
 \end{align*}

Consider the summand $ K_2(t,\tau), $ using the explicit
representation of the Green functions
\begin{align*}
 & K_2(t,\tau)=\frac1{4\pi}\imx\int\limits_0^{\tau}\frac{\mu_1'(\s1)-f(0,\s1)}
 {((\th_1(t)-\th_1(\tau))(\th_2(\tau)-\th_2(\s1)))^{3/2}}\\
 &\times
 \sum\limits_{n,m=-\infty}^{\infty}(\xi+2nh)(\xi+2mh)\exp\left(-\frac
 {(\xi+2nh)^2}{4(\th_1(t)-\th_1(\tau))}-\frac
 {(\xi+2mh)^2}{4(\th_2(\tau)-\th_2(\s1))}\right) d\s1 d\xi.
 \end{align*}
 Separating out of the series the term which corresponds to $ n=0 $
and $ m=0 $, we shall estimate the expression
\begin{align*}
  &K_{2,0}=\frac1{4\pi((\th_1(t)-\th_1(\tau))(\th_2(\tau)-\th_2(\s1)))^
  {3/2}} \times\\
 &\times
 \imx \xi^2\exp\left(-\frac{\xi^2}{4(\th_1(t)-\th_1(\tau))}-\frac{\xi^2}
 {4(\th_2(\tau)-\th_2(\s1))}\right)d\xi.
  \end{align*}

Changing the variable in the latter integral $
z=\dfrac{\xi}2\sqrt{\dfrac
{\th_1(t)-\th_1(\tau)+\th_2(\tau)-\th_2(\s1)}{(\th_1(t)-\th_1
(\tau))(\th_2 (\tau)-\th_2(\s1))}},$ we have
 $$
 K_{2,0}=\frac2{\pi(\th_1(t)-\th_1(\tau)+\th_2(\tau)-\th_2(\s1))^{3/2}}\int
 \limits_0^{T(t,\tau,\sigma)} z^2\exp(-z^2)dz,
 $$
 where $ T(t,\tau,\sigma)=\dfrac{h}2\sqrt{\dfrac{\th_1(t)-\th_1
 (\tau)+\th_2(\tau)-\th_2(\s1)}{(\th_1(t)-\th_1(\tau))(\th_2(\tau)-\th_2
 (\s1))}}$.
 
Integrating by parts we reduce the previous expression to the form
 \begin{align*}
 \adb
&K_{2,0}=-\frac{h}{2\pi(\th_1(t)-\th_1(\tau)+\th_2(\tau)-\th_2(\s1))\sqrt
 {(\th_1(t)-\th_1(\tau))(\th_2(\tau)-\th_2(\s1))}}\\
&\times\exp\left(-\frac{h^2}4\left
 (\frac1{\th_1(t)-\th_1(\tau)}+\frac1{\th_2(\tau)-\th_2(\s1)}\right)
 \right)+\\
&+ \frac1
 {\pi(\th_1(t)-\th_1(\tau)+\th_2(\tau)-\th_2(\s1))^{3/2}}\int\limits_0^
 {T(t,\tau,\s1)}\exp(-z^2)dz.
 \end{align*}
Estimating $ K_{2,0} $, we obtain the inequality
 $$
 |K_{2,0}|\le \frac{C_{16}}{(\th_1(t)-\th_1(\tau)+\th_2(\tau)-\th_2(\s1))^{3/2}}.
  $$
Applying (\ref{3.10}) we have
 $$
 |K_{2,0}|\le \frac{C_{17}}{(t^{\bb}-\s1^{\bb})^{3/2}}.
 $$
 Obviously, the estimates of other summands in $ K_2(t,\tau) $ are similar.
Returning to the estimate of $ K_2(t,\tau) $ and using the notation
(\ref{3.5}), we obtain
\begin{align*}
\adb
&\left|\frac{a_1(t)\tau^{\beta}}{t^{\beta}u_{2x}(0,t)}K_2(t,\tau)\right|\le
\frac{C_{18}A_1\tau^{\beta}}{\imi\frac{d\tau}{\m1}}\int
\limits_0^{\tau}\frac{d\s1}{(t^{\bb}-\s1^{\bb})^{3/2}}\\
&\le C_{19}
\tau^{\beta}\frac{t^{(\beta-1)/2}}{I_1}\int\limits_0^{\tau}\frac{d\s1}
{(t^{\bb}-\s1^{\bb})^{3/2}}.
\end{align*}
Consider the integral in the latter inequality:
\begin{align*}
&\int\limits_0^{\tau}\frac{d\s1}{(t^{\bb}-\s1^{\bb})^{3/2}}=\frac1{t^{3(\beta+1)
/2}}\int\limits_0^{\tau}\frac{d\s1}{\left(1-(\s1/t)^
{\bb}\right)^{3/2}}\le\frac1{t^{3\beta/2}}\int\limits_0^{\tau}\frac{d\s1}
{(t-\s1)^{3/2}}\le\\
&\le\frac2{t^{3\beta/2}\sqrt{t-\tau}}.
\end{align*}
 Finally we have
$$
\left|\frac{a_1(t)\tau^{\beta}}{t^{\beta}u_{2x}(0,t)}K_2(t,\tau)\right|
\le\frac{C_{20} \tau^{\beta}t^{(\beta-1)/2}}{t^{3\beta/2}\sqrt
{t-\tau}}\le\frac{C_{21}}{\sqrt{t(t-\tau)}}.
$$

Other summands $ K_{i}(t,\tau) $ are estimated analogously.  It
follows that  we get such inequality for $ K(t,\tau): $
$$
|K(t,\tau)|\le \frac{C_{22}}{\sqrt{t(t-\tau)}}.
$$
Taking into account  that $ a_0(t)\equiv0 $ for  $ t\in[0,\tilde
t],$ we finally have
$$
|K(t,\tau)|\le\frac{C_{22}}{\sqrt{\tilde t(t-\tau)}}\le\frac{C_{23}}
{\sqrt{t-\tau}},\quad t\in[0,T].
$$

This means that the equation (\ref{5.12}) as a homogeneous Volterra
integral equation of the second kind has only trivial solution $
a_0(t) \equiv 0.$ Then $ a(t) \equiv 0, t\in[0,T] $ and $
u(x,t)\equiv 0, (x,t)\in\overline Q_T $  as a consequence of the
uniqueness of the solution of direct problem
(\ref{5.5})-(\ref{5.7}).

 {\bf Remark 1.} The theorem of existence and uniqueness of solution for the
 problem (\ref{1.1})-(\ref{1.4}) may be expanded to the problem with another boundary and
 overdetermination conditions. Really, consider the analogous problem for
 equation (\ref{1.1}) with initial condition (\ref{1.2}), boundary conditions
 $ u_x(0,t)=\mu_1(t) ,\lb u_x(h,t)=\mu_2(t) $ and overdetermination condition
  $ u(0,t)=\mu_3(t)$. Then by change of unknown function  $ u_x\equiv v$
the given problem is reduced to the following one:
 \adb
 \begin{align*}
 & v_t=a(t)v_{xx}+f_x(x,t),\quad (x,t)\in Q_T,\\
 & v(x,0)=\f1'(x), \quad x\in[0,h],\\
& v(0,t)=\mu_1(t),\quad  v(h,t)=\mu_2(t), \quad t\in[0,T],\\
& a(t)v_x(0,t)=\mu_3'(t)-f(0,t), \quad t\in[0,T].
\end{align*}

 {\bf Remark 2.}  The assumptions of the theorem on functions $ \f1 $ and
 $ f $ may be relaxed and reduced to the following conditions: $ \f1\in
 C^1[0,h], f\in C^{1,0}(\overline Q_T).$ To check this statement, it is
 sufficient to study the behavior of the corresponding summands in
 (\ref{5.14}).

\pagebreak

\bigskip
{\bf References}
\bigskip

\begin{description}
\item[1.] Waid M C 1974 The first initial-boundary value problem for some
nonlinear time degenerate parabolic equation {\it Proc. Amer. Math.
Soc.} {\bf 42} 61-84
\item[2.]  Caffarelli L A, Friedman A 1979 Continuity of the density of a
gas flow in a porous medium {\it Trans. Amer. Math. Soc.} {\bf 252}
99-113
\item[3.] Podgayev A G 1987 On initial-boundary value problems for some
quasilinear  parabolic equations with nonclassical  degeneration
{\it Sibirsk. mat. journ. } {\bf 28} 129-139
\item[4.] Grebenev V.N. 1994 On a system of degenerate parabolic equations
that arises in fluid dynamics {\it Sibirsk. mat. journ. } {\bf 35}
753-767
\item[5.] Berestycki H, Busca J, Florent I 2000  An inverse  parabolic
problem arising in finance {\it C. R. Acad. Sci. Paris} {\bf 331}
965-969
\item[6.] Glushko V P  1968 Degenerate linear differential equations
{\it Diff. Uravneniya } {\bf 4 } 1957-1966
\item[7.] Kalashnikov A S 1968 On growing solutions of second order linear
equations with non-negative characteristic form {\it Mat. Zametki }
{\bf 3}  171-178.
\item[8.] Glushak A V, Shmulevich S D 1986 On some well-posed problems for higher
order parabolic equations which degenerates on time variable  {\it
Diff. Uravneniya } {\bf 22 } 1065-1068
\item[9.] Matsuzawa T. 1971 Sur une classe d'equations paraboliques degenerees
{\it Ann. sci. Ecole norm. super. } {\bf 104} 1-19
\item[10.] Malyshev I 1991 On the parabolic potentials in degenerate-type heat
equation {\it J. Appl. Math. Stochastic Anal} {\bf 4} 147-160
\item[11.] DiBenedetto E 1993 {\it Degenerate Parabolic Equations} (New York:
Springer)
\item[12.] Beland Y 2001 Time-vanishing properties of solutions of some degenerate
parabolic equations with strong absorption {\it Adv. Nonlinear
Stud.} {\bf 1} 117-152
\item[13.] Isakov V 1997 {\it Inverse Problems for Partial Differential Equations}
(New York: Springer)
\item[14.] Prilepko A I, Solovyev V V 1987 On solvability of inverse problems for
determination of lower order coefficient in a parabolic equation
{\it Diff. Uravneniya } {\bf 23 } 136-143
\item[15.] Pilant M, Rundell W 1988 Fixed point method for a nonlinear parabolic
inverse coefficient problem {\it Commun. Part. Dif. Equat. } {\bf
13} 469-493
\item[16.] Kirsch A 1996 {\it An Introduction to the Mathematical Theory of
Inverse Problems} (New York: Springer)
\item[17.]  Ivanchov M 2003 {\it Inverse Problems for Equations of Parabolic
Type } (Lviv: VNTL Publishers)
\item[18.] Eldesbaev T 1987 An inverse problem for a second order degenerate
hyperbolic equation {\it Izv. Acad. Nauk Resp. Kazakhstan. Ser. Fiz.
Mat.} {\bf 3} 27-29
\item[19.] Gadjiyev M M 1987 An inverse problem for degenerate elliptic equation
 {\it Primeneniya Metodov Funk. Anal. v Uravneniyah Mat. Fiz. } (Novosibirsk)
 66-71
\item[20.] Eldesbaev T 1995 On an inverse problem for a parabolic equation
whose order degenerates. I {\it Izv. Nats. Acad. Nauk Resp.
Kazakhstan. Ser.
 Fiz. Mat.} {\bf 3} 37-41
\item[21.] Ivanchov M, Saldina N An inverse problem for the degenerate heat
equation (to appear)
\item[22.] Jones B F 1962 The determination of a coefficient in a parabolic
differential equation. Part I. Existence and uniqueness {\it J.
Math. Mech.} {\bf 11} 907-918
\end{description}

\end{document}